\pdfoutput=1
\documentclass[a4paper,11pt]{article}
\usepackage[english]{babel}
\usepackage{amsfonts}
\usepackage{amsthm}
\usepackage{amsmath}
\usepackage{amssymb}
\usepackage{enumerate}
\usepackage{mathtools}
\usepackage{tikz}
\usepackage{stmaryrd}
\expandafter\def\csname opt@stmaryrd.sty\endcsname
{only,shortleftarrow,shortrightarrow}
\usepackage{tikz-cd}
\usepackage[all]{xy}
\usepackage{graphicx}
\usepackage{amsmath}
\usepackage{adjustbox}
\usepackage{graphicx}
\usepackage{extpfeil}
\usepackage{comment}
\usepackage{mathabx}
\usepackage{float}
\usepackage[labelformat=empty]{caption}
\usepackage[toc]{appendix}
\usepackage[left=3.2cm,top=3cm,right=3.2cm,bottom=3cm]{geometry}
\usepackage{hyperref}
\hypersetup{
colorlinks=true,
linkcolor=blue,
filecolor=magenta,      
urlcolor=cyan,
citecolor=[RGB]{0,204,0},
}
\usepackage{resizegather}
\usepackage[activate={true, nocompatibility}, final, tracking=true, kerning=true, spacing=true, factor=1100, stretch=10, shrink=10]{microtype}
\usepackage{dsfont}
\usetikzlibrary{matrix}
\usepackage[nameinlink]{cleveref}
\usepackage[normalem]{ulem}
\usepackage{sidecap}
\usepackage{mathrsfs}
\usepackage{mathscinet}

\usepackage{color} %CLASHES with hyperref
\definecolor{darkgreen}{rgb}{0.0, 0.7, 0.0}
\newenvironment{??}{\noindent \color{darkgreen}{\bf ???:} \footnotesize}{}
\definecolor{cyan}{cmyk}{1,0,0,0}

\newcommand{\bdg}{\begin{dg}}

\theoremstyle{plain}
\newtheorem{thm}{Theorem}
\newtheorem{coroll}[thm]{Corollary}

\newtheorem{prop}[thm]{Proposition}

\newtheorem{remark}[thm]{Remark}

\newtheorem*{thm*}{Main Theorem}

\theoremstyle{definition}

\tikzset{
  symbol/.style={
    draw=none,
    every to/.append style={
      edge node={node [sloped, allow upside down, auto=false]{$#1$}}}
  }
}
\microtypecontext{spacing=nonfrench}

\newcommand{\coker}{\mathrm{coker}}
\newcommand{\id}{\mathrm{id}}

\newcommand{\cO}{\mathcal{O}}

\newcommand{\cE}{\mathcal{E}}
\newcommand{\cF}{\mathcal{F}}
\newcommand{\cG}{\mathcal{G}}

\newcommand{\cK}{\mathcal{K}}
\newcommand{\cH}{\mathcal{H}}

\newcommand{\cQ}{\mathcal{Q}}
\newcommand{\cV}{\mathcal{V}}

\newcommand{\bA}{\mathbb{A}}

\newcommand{\on}{\operatorname}

\renewcommand{\ker}{ \on{ker}}

\newcommand{\Spec}{\on{Spec}}

\newcommand{\GL}{\mathrm{GL}}

\usepackage{color} %CLASHES with hyperref
\definecolor{darkgreen}{rgb}{0.0, 0.7, 0.0}

\begin{document}
\title{\textbf{Relative cohomology and base-change}}
\author{Andres Fernandez Herrero and Giovanni Inchiostro}
\date{}

\maketitle
\begin{abstract}
    In this note we prove a simple generalization of cohomology and base-change which seems to be missing from standard sources. We illustrate its usefulness by providing an alternative proof of a result of Ishimura \cite{ishimura_descent} on descent of vector bundles along blowups at a smooth center.
\end{abstract}

\section{Introduction}
In this note we adapt arguments from \cite[\S 7]{egaiii_2} to prove a relative version of cohomology and base-change, see \Cref{thm: main theorem on pushforward and base-change} and \Cref{prop: vanishing of higher pushfowards fiberwise}. Our main motivation is to have a citable source for these results, which we have found to be useful in various contexts such as:
\begin{itemize}
    \item Computing pushforwards of vector bundles by inducting on the dimension. We illustrate this strategy by providing an alternative proof of a result of Ishimura \cite{ishimura_descent} on descent of vector bundles along blowups at a smooth center (see \Cref{prop: application descent vector bundes blowups}).
    \item Showing that certain constructions in moduli theory can be performed at the level of families. This has been used in several projects by the authors \cite{gauged-maps-projective, alternative_moduli_curves, gieseker_sl_n_bundles, fundamental_groups_vector_bundles}. We also record an application to Stein factorizations in families in this note (see \Cref{prop: stein factorization in families}) 
\end{itemize}

Our results in this note are straightforward generalizations of the usual statement of cohomology and base-change, variations of which can be found in many sources \cite{egaiii_2, hartshorne_book, lange, stacks-project, gortz-wedhorn-ii, vakil_rising_sea} .

\subsubsection*{Assumptions and notation.}
Unless otherwise stated, we work over a fixed locally Noetherian base scheme $S$. Given two morphisms of schemes $X \to S$ and $T \to S$, we may sometimes use the notation $X_T: = X\times_S T$ to the denote the base-change. In this case, given a quasicoherent sheaf $\cF$ on $X$, we denote by $\cF_T$ its pullback to $X_T$.

\subsection*{Acknowledgements.} We thank Francisco Gallardo for pointing out a mistake in a previous formulation of one of the results. We also thank Nicol\'as Vilches Reyes and Jean-Marc Drezet for pointing out relevant references in the literature.

\section{Main results}

\subsection{A result on derived pushfowards}\label{subsection_result_on_derived_pushfowards}
\begin{thm} \label{thm: main theorem on pushforward and base-change}
    Let $f: X \to Y$ be a finite type separated morphism of schemes which are locally of finite type over $S$. Let $\cF^{\bullet}$ be a bounded complex of quasicoherent sheaves on $X$ such that $\cF^i$ is $S$-flat for all $i$. Then, the following hold:
    \begin{enumerate}[(a)]
        \item The derived pushforward $Rf_*(\cF^{\bullet})$ is Zariski locally on $Y$ quasi-isomorphic to a bounded complex of $S$-flat quasicoherent $\cO_Y$-sheaves.
        \item Let $g:T \to S$ be a morphism of locally Noetherian schemes, and consider the induced Cartesian diagram
             \[
\begin{tikzcd}
 X_T \ar[r, "g_X"]
  \ar[d, "f_T"] & X \ar[d, "f"] \\ Y_T \ar[r, "g_Y"] &
  Y.
\end{tikzcd}
\]
If $Y \to S$ is flat, then there is a natural quasi-isomorphism $Lg_Y^*(Rf_*(\cF^{\bullet})) \to R(f_T)_*(g_X^* \cF^{\bullet})$, where $g_X^*(-)$ denotes the naive pullback of the complex (as opposed to the derived pullback).
\item If the morphism $Y \to S$ is flat and for all integers $i$ we have that $\cH^i(\cF^{\bullet})$ is coherent with proper support over $Y$, then $Rf_*(\cF^{\bullet})$ is Zariski locally on $Y$ quasi-isomorphic to a bounded complex of $S$-flat coherent $\cO_Y$-sheaves.
    \end{enumerate}
\end{thm}
\begin{remark}
    Note that, contrary to the usual cohomology and base-change setting, we are not requiring $f$ to be flat.
\end{remark}
\begin{proof}
    All the statements can be checked Zariski-locally on $Y$, and hence we may assume without loss of generality that both $S$ and $Y$ are affine. Say $Y = \Spec(B)$ and $S = \Spec(A)$. Choose a finite cover $X= \bigcup_{j \in J} X_j$ by affine open subschemes $X_j \subset X$. For every $\Sigma \subset J$, set $X_{\Sigma} = \bigcap_{j \in \Sigma} X_j$. Then the derived pushforward $Rf_*(\cF^{\bullet})$ is naturally quasi-isomorphic to the totalization $Tot(\cG^{\bullet, \bullet})$ of the double complex $\cG^{n,m} := \oplus_{\Sigma \subset J, |\Sigma| = n} \cF^m(X_{\Sigma})$ of $B$-modules. By the assumption on $\cF^{\bullet}$, every term in the double-complex $\cG^{n,m}$ is $A$-flat. This shows part (a).

    For part (b), there is a natural base-change morphism $Lg_Y^*(Rf_*(\cF^{\bullet})) \to R(f_T)_*(Lg_X^* \cF^{\bullet})$ induced by adjunction. Composing with the morphism $Lg_X^*\cF^{\bullet} \to g_X^*\cF^{\bullet}$, we obtain the required natural morphism $Lg_Y^*(Rf_*(\cF^{\bullet})) \to R(f_T)_*(g_X^* \cF^{\bullet})$. To check that this is a quasi-isomorphism, we may work locally on $Y$ and assume that we have $Y = \Spec(B)$, $S = \Spec(A)$ and $T= \Spec(C)$. Note that the base-change $X_T = \bigcup_{j \in J} X_j \times_S T$ is again a finite affine cover of the base-change $X_T$. By the same reasoning as in part (a), the complex $R(f_T)_*(g_X^* \cF^{\bullet})$ is naturally quasi-isomorphic to the totalization of the double complex of $B \otimes_A C$-modules $\cG_T^{n,m} := \oplus_{\Sigma \subset J, |\Sigma| = n} g_X^*\cF^m(X_{\Sigma} \times_S T) = \oplus_{\Sigma \subset J, |\Sigma| = n} \cF^m(X_{\Sigma}) \otimes_A C$.

    Since every term in the complex $Tot(\cG^{\bullet, \bullet})$ representing $Rf_*(\cF^{\bullet})$ is $S$-flat and $Y \to S$ is flat, we have that $Lg_Y^*(Rf_*(\cF^{\bullet})$ is naturally quasi-isomorphic to the naive pulled back complex $Tot(\cG^{\bullet, \bullet})\otimes_A C$ of $B \otimes_A C$-modules. Under these identifications the morphism $Lg_Y^*(Rf_*(\cF^{\bullet})) \to R(f_T)_*(Lg_X^* \cF^{\bullet})$ is seen to induce the natural isomorphism of complexes $Tot(\cG^{\bullet, \bullet}\otimes_A C) \xrightarrow{\sim} Tot(\cG^{\bullet}) \otimes_A C$. This concludes the proof of (b).

    For the proof of (c), we note that, as by assumption $\cH^i(\cF)$ has proper support over $Y$, the cohomologies of the bounded complex $Tot(\cG^{\bullet, \bullet})$ of $A$-flat $B$-modules are finitely generated $B$-modules. Now the desired statement follows from \cite[Lemma 1 page 47]{Mum74} and its proof. Indeed, from \cite[Lemma 1 page 47]{Mum74} the complex $Tot(\cG^{\bullet, \bullet})$ is quasi-isomorphic to a bounded complex of finitely generated $B$-modules of the form \[\cV^\bullet\colon\text{ }\text{ } 0\to \cV_d\to \cV_{d+1}\to \ldots \to \cV_{d'}\to 0\] for a certain $d$, such that $\cV_j$ are free for $j>d$. One can then argue that $\cV_d$ is $A$-flat exactly as in \cite[Lemma 1 page 47]{Mum74}, by explicitly constructing the mapping cylinder of $\cV^\bullet \to Tot(\cG^{\bullet, \bullet})$ and using the $A$-flatness of the free $B$-modules $\cV_j$ for $j>d$.
\end{proof}

\subsection{Relative cohomology and base-change}
In this subsection we record some standard consequences of \Cref{thm: main theorem on pushforward and base-change}. The relevant statements in \cite[\S 7]{egaiii_2} do not quite apply on the nose, so we explain the minor modifications needed to carry out the proofs in the relative setting. 

\begin{prop} \label{prop: vanishing of higher pushfowards fiberwise}
    Let $f: X \to Y$ be a separated finite type morphism of schemes which are locally of finite type over $S$. Assume that $Y \to S$ is flat. Let $\cF$ be a coherent $S$-flat sheaf on $X$ whose support is proper over $Y$. Let $s \in S$, and consider the morphism of fibers $f_s: X_s \to Y_s$. Suppose that we have $R^i(f_s)_*(\cF_s) =0$ for some $i\geq 0$. Then, the following hold
    \begin{enumerate}[(a)]
        \item There exists some open subscheme $U \subset Y$ containing the fiber $Y_s$ such that if we set $W:= f^{-1}(U)$, then the morphism $f: W \to U$ satisfies that for all morphisms $T \to S$ we have $R^i(f_T)_*(\cF|_{W_T}) =0$.
        \item With notation as in part (a), the formation of the coherent sheaf $R^{i-1}f_*(\cF|_{W})$ on $U$ commutes with arbitrary $S$ base-change, and the following are equivalent:
        \begin{itemize}
            \item $R^{i-1}f_*(\cF|_{W})$ is $S$-flat.
            \item The formation of the coherent sheaf $R^{i-2}f_*(\cF|_{W})$ on $U$ commutes with arbitrary $S$ base-change.
        \end{itemize}
    \end{enumerate}
\end{prop}
\begin{proof}
    Since everything can be checked locally on $Y$, we may assume that both $Y$ and $S$ are affine, say $Y=\Spec(B)$ and $S= \Spec(A)$. In view of \Cref{thm: main theorem on pushforward and base-change}(c), we may represent $Rf_*(\cF)$ by a bounded complex $\cE$ of finitely generated $B$-modules that are $A$-flat. Furthermore, in view of \Cref{thm: main theorem on pushforward and base-change}(b), for all $T=\Spec(C) \to S$ we have that $R(f_T)_*(\cF_T)$ is quasi-isomorphic to the complex $\cE^{\bullet} \otimes_A C$ obtained by naively tensoring with $C$. 

Let us prove part (a). Let $\overline{\cE}^i := \coker \left( \partial: \cE^{i-1} \to \cE^i\right)$ be the cokernel of the differential. The next differential $\partial: \cE^i \to \cE^{i+1}$ induces a morphism $\overline{\partial}: \overline{\cE}^i \to \cE^{i+1}$. For all $A$-algebras $C$, we have 
\[ \overline{\cE}^i \otimes_A C = \coker \left( \partial \otimes_A id_C : \cE^{i-1} \otimes_A C\to \cE^i\otimes_A C\right), \; \; \text{and}\]
\[H^i(\cE^{\bullet} \otimes_A C) = \ker\left(\overline{\partial}\otimes_A id_C: \overline{\cE}^i\otimes_A C \to \cE^{i+1} \otimes_A C\right).\] 
Let $A_s$ denote the local ring of $A$ at $s$, and let $k(s)$ denote the residue field of $s$. Our assumption $R^i(f_s)_*(\cF_s) =0$ translates into $\ker(\overline{\partial}\otimes_A id_{k(s)}) =0$. Since $\cE^{i+1}$ is $A$-flat, it follows from the slicing criterion for flatness \cite[\href{https://stacks.math.columbia.edu/tag/00ME}{Tag 00ME}]{stacks-project} that the localized morphism $\overline{\partial}\otimes_A A_s: \overline{\cE}^i \otimes_A A_s \to \cE^{i+1} \otimes_A A_s$ is injective with $A_s$-flat cokernel. Since $\overline{\cE}^i$ is a finitely generated $B$-module and $B$ is Noetherian, there is a single element $f \in A$ such that $s \in \Spec(A_f) \subset \Spec(A)$ and the localization $\overline{\partial}: \overline{\cE}^i \otimes_A A_f \to \cE^{i+1} \otimes_A A_f$ is injective. Replacing $S = \Spec(A)$ with $\Spec(A_f)$, we may assume that $\overline{\partial}$ is injective. Consider the cokernel $\cQ^i := \coker(\overline{\partial})$. We have shown that $\cQ^i \otimes_A A_s$ is $A_s$-flat, and hence the sheaf $\cQ^i$ on $Y= \Spec(B)$ is $A$-flat at every point of the fiber $Y_s$. Since the $A$-flat locus of $\cQ^i$ is open \cite[\href{https://stacks.math.columbia.edu/tag/0399}{Tag 0399}]{stacks-project}, after replacing $Y=\Spec(B)$ with an open subscheme containing the fiber $Y_s$, we may assume without loss of generality that the injective morphism $\overline{\partial}: \overline{\cE}^i \to \cE^{i+1}$ has $A$-flat cokernel. It follows that the morphism remains injective after arbitrary base-change, and hence for all $A$-algebras $C$ we have 
\[H^i(\cE^{\bullet} \otimes_A C) = \ker\left(\overline{\partial}\otimes_A id_C: \overline{\cE}^i\otimes_A C \to \cE^{i+1} \otimes_A C\right)=0.\]
The statement in part (a) follows, since the complex $\cE^{\bullet}$ universally computes the higher direct images of $f$.

We are left to show part (b). For ease of notation, we replace $Y$ with the open $W \subset Y$ we have found in part (a). Using the notation above, the $A$-flatness of $\cQ^i$ shows that $\overline{\cE}^i= \coker(\partial: \cE^{i-1} \to \cE^i)$ is $A$-flat, and that for all $A$-algebras $C$ the morphism $\overline{\partial}\otimes_A id_C: \overline{\cE}^i \otimes_A C \to \cE^{i+1} \otimes_A C$ is injective. Hence, if we set $\cK^i$ to be the kernel of the surjection $\cE^i \twoheadrightarrow \overline{\cE}^i$, then $\cK^i$ is $A$-flat, and we have 
\begin{align*}
    \cK^i \otimes_A C & = \ker\left(\cE^i\otimes_A C \to \overline{\cE}^i \otimes_A C\right)\\
    & =\ker\left(\partial\otimes_A \id_C: \cE^i\otimes_A C \to \overline{\cE}^i \otimes_A C \hookrightarrow \cE^{i+1} \otimes_A C\right).
\end{align*}
for any $A$-algebra $C$. This shows that the cohomological truncation
\[ \tau_{\leq i}(\cE^{\bullet}) = \left[ \ldots \to \cE^0 \to \cE^1 \to \ldots \to \cE^{i-1} \to \cK^i \to 0 \to \ldots\right]\]
is a complex with $A$-flat terms, and that its formation commutes with $A$ base-change. Hence the complex $\tau_{\leq i}(\cE^{\bullet})$ universally computes the higher direct images $R^jf_*(\cF)$ for $j \leq i$. Furthermore, we know that $\cE^{i-1} \to \cK^i$ is surjective, since $H^i(\tau_{\leq i}(\cE^i)) = R^if_*(\cF) =0$. Hence by the same reasoning, if we set $\cK^{i-1}$ to be the kernel of $\cE^{i-1} \to \cK^i$, then the cohomological truncation
\[ \tau_{\leq i-1}(\cE^{\bullet}) = \left[ \ldots \to \cE^0 \to \cE^1 \to \ldots \to \cE^{i-2} \to \cK^{i-1} \to 0 \to \ldots\right]\]
has $A$-flat terms and its formation commutes with $A$ base-change. We have that $R^{i-1}f_*(\cF) = \coker(\cE^{i-2} \to \cK^{i-1})$, and its formation commutes with $A$ base-change by the right exactness of tensoring. The equivalence of the bullet points in part (b) is an immediate consequence of the equivalence of the first two items in \cite[Prop. 7.4.2]{egaiii_2} applied to the complex of $A$-modules $\tau_{\leq i-1}(\cE^{\bullet})$.
\end{proof}

\begin{coroll} \label{coroll: vanishing of first pushforward in fibers and flatness}
    Let $f: X \to Y$ be a separated finite type morphism of schemes which are locally of finite type over $S$. Assume that $Y \to S$ is flat. Let $\cF$ be a coherent $S$-flat sheaf on $X$ whose support is proper over $Y$. Suppose that for all $s \in S$, we have $R^1(f_s)_*(\cF_s)=0$. Then, the following hold
    \begin{enumerate}[(a)]
        \item The pushforward $f_*(\cF)$ is an $S$-flat coherent sheaf on $Y$, and its formation commutes with arbitrary $S$ base-change.
        \item $f_*(\cF)$ is locally free if and only if for all $s \in S$ we have that $(f_s)_*(\cF)$ is a locally free sheaf on $Y_s$.
    \end{enumerate} 
\end{coroll}
\begin{proof}
    Part (a) follows directly from \Cref{prop: vanishing of higher pushfowards fiberwise}. Part (b) follows from part (a) and the fiberwise criterion for flatness \cite[\href{https://stacks.math.columbia.edu/tag/039C}{Tag 039C}]{stacks-project} applied to the morphism $Y \to S$.
\end{proof}

% \begin{prop} \label{prop: fiberwise local constancy of fibers}
%     Let $f: X \to Y$ be a separated finite type morphism of schemes which are locally of finite type over $S$. Assume that $Y \to S$ is flat, and $Y$ is reduced. Let $\cF$ be a coherent $S$-flat sheaf on $X$ whose support is proper over $Y$. Fix $i \geq 0$. Suppose that there exists some integer $r$ such that for all $s \in S$ we have that $R^i(f_s)_*(\cF_s)$ is a locally free sheaf of rank $r$ on $Y_s$. Then, the following hold:
%     \begin{enumerate}[(1)]
%         \item The formation of both $R^{i}f_*(\cF)$ and $R^{i-1}f_*(\cF)$ commutes with arbitrary $S$ base-change.
%         \item $R^if_*(\cF)$ is a locally free sheaf of rank $r$ on $Y$.
%     \end{enumerate}
% \end{prop}
% \begin{proof}
%     \andres{TO DO: We don't really need this, but I thought that it could potentially be useful to add this, since it is a statement that is usually included in the package of cohomology and base-change. Proof to be added.}
% \end{proof}

\section{Some applications}

\begin{prop} \label{prop: stein factorization in families}
    Let $f: X \to Y$ be a proper morphism of finite type flat $S$-schemes. If for all $s \in S$ we have $R^1(f_s)_*\cO_{X_s}=0$, then the Stein factorization $\underline{\Spec}_{Y}(f_*(\cO_X)$ is $S$-flat and its formation commutes with arbitrary $S$ base-change.
\end{prop}
\begin{proof}
    This is an immediate consequence of \Cref{coroll: vanishing of first pushforward in fibers and flatness}.
\end{proof}

Another way in which relative cohomology and base-change could be useful is in computing pushforwards of vector bundles under proper (not necessarily flat) morphisms $X \to Y$ of varieties by locally choosing a morphism $Y \to S$ to a smooth curve $S$ and looking at the fibers $X_s \to Y_s$. This allows for inductive arguments where at each step we decrease the dimension of the varieties involved.
We illustrate this kind of computation by providing a more elementary proof of \cite[Thm. 1]{ishimura_descent}. We remark that our proof is not significantly easier from the conceptual point of view, but it has the advantage of using only basic sheaf cohomology, thus avoiding the theorem on formal functions.
\begin{prop} \label{prop: application descent vector bundes blowups}
    Let $Y$ be a smooth scheme over a field $k$, and let $Z \subset X$ be a smooth closed subscheme. Let $X:=\operatorname{Bl}_Z(X)$ be the blowing up of $Z$ in $Y$ with $f\colon X \to Y$ the blow-up morphism, and let $E:= f^{-1}(Z)$ denote the exceptional divisor. Then, for all vector bundles $\cF$ on $X$, the following are equivalent:
    \begin{enumerate}[(a)]
        \item $\cF$ is isomorphic to the pullback on $Y$.
        \item The restriction $\cF|_{E}$ is isomorphic to the pullback of a vector bundle on $Z$.
    \end{enumerate}
\end{prop}
Note that this also follows as such a blow-up is a pushout in algebraic stacks (so in particular it is a pushout for maps to $B \GL_n$) from \cite[Theorem 1.1]{arena2023criterion}.
\begin{proof}
    The implication (a) $\Rightarrow$ (b) is clear. Let us show the converse. Assume that $\cF|_{E}$ is the pullback of a vector bundle on $Z$. We shall prove the following stronger property of the pair $(\cF, f)$:
    
    \textbf{(*)} We have $R^1f_*(\cF)=0$, the pushforward $f_*(\cF)$ is a vector bundle, and the counit morphism $f^*f_*(\cF) \to \cF$ is an isomorphism.

    These properties can be checked \'etale locally on $Y$, and hence we may assume without loss of generality that $\cF|_{E}$ is trivializable, that $Y$ is an open subset of an affine space $\mathbb{A}^n$ with coordinates $x_1, x_2, \ldots , x_n$, and that $Z \subset Y$ is given by the vanishing of the coordinates $x_1, x_2, \ldots, x_m$ for some $m \le n$. We prove \textbf{(*)} by induction on the dimension $n$. The base-case of dimension $n=1$ is trivial.
    
    Consider the morphism $Y \hookrightarrow \mathbb{A}^n \to \mathbb{A}^{n-m}$ obtained by projection to the last $n-m$-coordinates. Observe that, from \cite[\href{https://stacks.math.columbia.edu/tag/00R4}{Tag 00R4}]{stacks-project}, the composition $X \xrightarrow{f} Y \to \mathbb{A}^{n-m}$ is flat. In view of \Cref{coroll: vanishing of first pushforward in fibers and flatness}, to conclude the proof of the claim, it suffices to show that for all $s \in \mathbb{A}^{n-m}$ the pair $(\cF_s, f_s: X_s \to Y_s)$ satisfies the property (*) above. But, for a given $s \in \mathbb{A}^{n-m}$, the fiber $Y_s$ is an open subset of an affine space of dimension $m$, and $f_s: X_s \to Y_s$ is the blowing up of the origin. Hence, we have reduced to the case when $Z$ is the origin in $\bA^n$. 
    
    In this case, we can further consider the morphism $Y \to \mathbb{A}^1$ obtained by projecting to the last coordinate. The composition $X \to Y \to \mathbb{A}^1$ is flat, and it suffices to show that for all $s \in \mathbb{A}^1$ the property \textbf{(*)} holds for the pair $(\cF_s, f_s: X_s \to Y_s)$.

    There are two cases. If $s \in \mathbb{A}^1$ is not the origin, then $f_s$ is an isomorphism and property (*) is clear. Otherwise, if $s$ is the origin, then $Y_s$ is an open subset of $\mathbb{A}^{n-1}$ and $X_s$ has two irreducible components $X_s = \Sigma_1 \cup \Sigma_2$, where $\Sigma_1$ is the blowing up of $Y_s$ at the origin of $\mathbb{A}^{n-1}$ and $\Sigma_2 \cong \mathbb{P}^{n-2}$. If we denote by $E' \subset \Sigma_1$ the exceptional divisor of the blowing up $\Sigma_1 \to Y$, then $\Sigma_2= \mathbb{P}^{n-2}$ is glued to $\Sigma_1$ via an isomorphism of the exceptional divisor $E' \cong \mathbb{P}^{n-3}$ with a hyperplane of $\mathbb{P}^{n-2}$. 

    We have a short exact sequence of sheaves on $X_s$
    \[ 0 \to \cO_{X_s} \to \cO_{\Sigma_1} \oplus \cO_{\Sigma_2} \to \cO_{E'} \to 0.\]
    Tensoring by $\cF_s$, we obtain a short exact sequence of sheaves
    \[ 0 \to \cF_s \to \cF_s|_{\Sigma_1} \oplus \cO^{\oplus r}_{\Sigma_2} \to \cO^{\oplus r}_{E'} \to 0,\]
    where we denote by $r$ the rank of $\cF$ and we are using that $\cF_s|_{\Sigma_2}$ is trivializable.
    Pushing forward via $f_s: X_s \to Y_s$ and using the evident fact that $\cO_{0} = R(f_s)_* \cO_{\Sigma_2} \to R(f_s)_*\cO_{E'} = \cO_0$ is an isomorphism of sheaves on $Y_s$ supported at the origin $0 \in Y_s$, we get a quasi-isomorphism $R(f_s)_*(\cF_s) \xrightarrow R(f_s)_*(\cF_s|_{\Sigma_1})$. In other words, $R(f_s)_*(\cF_s)$ is quasi-isomorphic to the derived pushforward of the vector bundle $\cF_s|_{\Sigma_1}$ via the blowing up morphism $\Sigma_1 \to Y_s$.
    
    Now, since $Y_s$ is an open inside of affine space of dimension $n-1$, we conclude that property (*) holds for the pair $(\cF_s|_{\Sigma_1}, \Sigma_1 \to Y_s)$ by induction. It follows that $R^1(f_s)_*(\cF_s)=0$ and $(f_s)_*(\cF_s)$ is a vector bundle. We get a corresponding commutative diagram of exact sequences of sheaves on $X_s$, where the vertical morphisms are induced by the counits
                 \[
\begin{tikzcd}
 0 \ar[r] & f_s^*(f_s)_*(\cF_s) \ar[r] \ar[d] & 
  \ar[d] \ar[r] f_s^*(f_s)_*(\cF_s|_{\Sigma_1}) \oplus f_s^*(\cO_0^{\oplus r}) & f_s^*(\cO_0^{\oplus r}) \ar[d] \\ 0 \ar[r] & \cF_s \ar[r] &  \cF_s|_{\Sigma_1} \oplus \cO^{\oplus r}_{\Sigma_2}  \ar[r] & \cO^{\oplus r}_{E'} \ar[r] & 0.
\end{tikzcd}
\]
Now, we established that the morphism $f_s^*(f_s)_*(\cF_s|_{\Sigma_1}) \to \cF_s|_{\Sigma_1}$ is an isomorphism by induction, and the morphism $f_s^*(\cO_0^{\oplus r}) \to \cO^{\oplus r}_{\Sigma_2}$ is an isomorphism by inspection. Since the morphism $f_s^*(f_s)_*(\cF_s|_{\Sigma_1}) \oplus f_s^*(\cO_0^{\oplus r}) \to f_s^*(\cO_0^{\oplus r})$ is evidently surjective, an application of the snake lemma implies that the counit $f_s^*(f_s)_*(\cF_s) \to \cF_s$ is an isomorphism, thus concluding the proof of property \textbf{(*)}.
\end{proof}

\footnotesize{\bibliography{moduli-sheaves.bib}}
\bibliographystyle{alpha}
\end{document}